\documentclass[10pt,twoside,reqno]{amsart}
\usepackage{mathrsfs}
\calclayout

\numberwithin{equation}{section}
\makeatletter

\renewcommand{\@secnumfont}{\bfseries}

\renewcommand{\section}{\@startsection{section}{1}%
  {0mm}{.7\linespacing\@plus\linespacing}{.5\linespacing}
  {\normalfont\bfseries\centering}}

\newcommand{\bibsection}{\@startsection{section}{1}%
  {0mm}{.7\linespacing\@plus\linespacing}{.5\linespacing}
  {\normalfont\scshape\centering}}

\renewcommand{\@biblabel}[1]{#1.}

\newtheorem{thm}{\bf Theorem}[section]

\begin{document}

\vspace{1.3cm}

\title
        {Identities of symmetry for higher-order generalized $q-$Euler polynomials}

\author{D. V. Dolgy, D. S. Kim, T. G. Kim, J. J. Seo }

\thanks{\scriptsize }
\address{1\\ Institute of Mathematics and Computer Sciences\\
             Far Eastern Federal University\\
              Vladivostok, 690060, Russia}
              \email{$d_{-}dol@mail.ru$}
\address{2\\ Department of Mathematics\\
             Sogang University\\
              Seoul 121-742, Republic of Korea}
\email{dskim@sogang.ac.kr}
\address{3\\ Jangjeon Research Institute for Mathematics and Physics\\
            Hapcheon, Gyoungnam, 678-800, Republic of Korea}
\email{tgkim2013@hotmail.com}
\address{4\\Department of Applied Mathematics\\
            Pukyong National University\\
            Busan 608-737, Republic of Korea}
\email{seo2011@pknu.ac.kr}

\keywords{higher-order generalited $q-$Euler polynomials}

\maketitle

\begin{abstract} In this paper, we investigate the properties of symmetry in two variables related to multiple Euler $q-l-$function which interpolates higher-order $q-$Euler polynomials at negative integers. From our investigation, we can derive many interesting identities of symmetry in two variables related to generalized higher-order $q-$Euler polynomials and alternating generalized $q-$power sums.
\end{abstract}

\pagestyle{myheadings}
\markboth{\centerline{\scriptsize D. V. Dolgy,  D. S. Kim, T. G.  Kim, J. J. Seo}}
          {\centerline{\scriptsize Identities of symmetry for higher-order generalized $q-$Euler polynomials}}
\bigskip
\bigskip
\medskip
\section{\bf Introduction}
\bigskip
\bigskip
Let $\chi$ be a Dirichlet character with $d\in{\mathbb N}$ with conductor $d\equiv1(mod2)$. Then the generalized Euler polynomials attached to $\chi$ are defined by the generating function to be\\
\begin{equation}
\ 2 \sum_{a=0}^{d-1}\frac{\chi(a)(-1)^a{e}^{(a+x)t}}{e^{dt}+1}=\sum_{n=0}^{\infty} E_{n,\chi}(x)\frac{t^n}{n!}, \ \ \textnormal{(see [4], [9], [19])}.
\end{equation}
The generalized Euler polynomials of order $r\in{\mathbb N}$ attached to $\chi$ are also defined by the generating function to be\\
\begin{equation}\begin{split}
\left(2\sum_{a=0}^{d-1}\frac{\chi(a)(-1)^a{e}^{(a+x)t}}{e^{dt}+1}\right)^r= \sum_{n=0}^{\infty} E_{n,\chi}^{(r)}(x)\frac{t^n}{n!}.
\end{split}\end{equation}
When $x=0, E_{n,\chi}^{(r)}=E_{n\chi}^{(r)}(0)$ are called the generalized Euler numbers attached to $\chi$, \ \ \textnormal{(see [9], [13])}.\\
Assume that $q\in{\mathbb C}$ with $|q|<1$ and define $q-$numbers by\\
\begin{equation}\begin{split}
\ [x]_q=\frac{1-q^x}{1-q},  \ \ \textnormal{(see [8]-[21])}.
\end{split}\end{equation}
Note the $\lim_{q\rightarrow1}[x]_q = x$.\\
In[12] and [13], Kim for the first time considered various $q$-extensions (or $(h,q)-$extensions) of Euler numbers and polynomials and constructed analytic continuations which interpolate his $q-$numbers and polynomials. Until recently, many authors have studied $q-$Euler or $(h,q)-$Euler polynomials due to him, \ \ \textnormal{(see [1]-[21])}. In [13], Kim defined the $(h,q)-$extension of generalized higher-order Euler polynomials attached to $\chi$ which is given by the generating function to be
\begin{equation}\begin{split}
\ F_{q,\chi}^{(h,r)}(t,x)&=[2]_q^r\sum_{m_1,\cdot\cdot\cdot,m_r=0}^{\infty} q^{\sum_{j=1}^{r}{(h-j+1)m_j}} {(-1)}^{\sum_{j=1}^{r}m_j}\left(\prod_{j=1}^{r}\chi(m_j)\right)e^{[x+{\sum_{l=1}^{r}m_l}]_q t}\\
&=\sum_{n=0}^{\infty} E_{n,\chi,q}^{(h,r)}(x)\frac{t^n}{n!},
\end{split}\end{equation}
where $ h\in {\mathbb Z}$ and $ r\in {\mathbb N}.$ \\
Note that\\
\begin{equation}\begin{split}
\lim_{q\rightarrow 1} F_q^{(h,r)}(t,x) &=\left(2\sum_{a=0}^{d-1}\frac{\chi(a)(-1)^a e^{(a+x)t}}{e^{dt}+1}\right)^r\nonumber\\
&=\sum_{n=0}^{\infty} E_{n,\chi}^{(r)}(x)\frac{t^n}{n!}.
\end{split}\end{equation}
When $x=0, E_{n,\chi,q}^{(h,r)} =E_{n,\chi,q}^{(h,r)}(0)$ are called the $(h,q)-$extension of generalized higher-order Euler numbers attached to $\chi$.\\
From(1.4), We note that\\
\begin{equation}\begin{split}
\ E_{n,\chi,q}^{(h,r)}(x) &=\sum_{l=0}^{n}{n\choose l} q^{lx} E_{l,\chi,q}^{(h,r)}[x]_q^{n-l}\\
 &=\left(q^x E_{\chi,q}^{(h,r)}+[x]_q\right)^n,
\end{split}\end{equation}
with the usual convention about replacing $\left(E_{\chi,q}^{(h.r)}\right)^n$ by $E_{n,\chi,q}^{(h.r)}$.\\
In[13], Dirichlet-type multiple $(h,q)-l-$ function is defined by Kim to be\\
\begin{equation}\begin{split}
\ l_{q,r}^{(h)}(s,x|\chi)=& \frac{1}{\Gamma(s)}\int_{0}^{\infty} F_{q,{\chi}}^{(h,r)}(-t,x)t^{s-1}dt\\
=&[2]_q^r\sum_{m_1,\cdot\cdot\cdot,m_r=0}^{\infty}\frac{q^{\sum_{l=1}^{r}(h-l+1)m_l} (\prod_{l=1}^r \chi(m_l)) (-1)^{\sum_{l=1}^r m_l}}{[m_1+\cdot\cdot\cdot+m_r+x]_q^s},
\end{split}\end{equation}
where ${s,h} \in{\mathbb C}$ and $ x\in{\mathbb R,}$ with $x\neq{0,-1,-2,\cdot\cdot\cdot}$.\\
By using Cauchy residue theorem, we get\\
\begin{equation}\begin{split}
\ l_{q,r}^{(h)}(-n,x|\chi)=  E_{n,{\chi},q}^{(h,r)}(x), \ \ \ n \in {\mathbb Z_{\geq0}}.
\end{split}\end{equation}
In this paper, we investigate the properties of symmetry in two variables related to Dirichlet-type multiple $(h,q)-$function which interpolates the $(h,q)-$extension of generalized higher-order Euler polynomials attached to $\chi$ at negative integers. From our investigation, we can derive many interesting identities of symmetry in two variables related to $(h,q)-$extension of generalized higher-order Euler polynomials and alternating generalized $q-$power sums.\\

\section{\bf Identities for the $(h,q)-$extension of generalized higher-order Euler polynomials }
\bigskip
\medskip
In this section, we assume that $\chi$ is a Dirichlet character with conductor $d\in {\mathbb N}$ with $d\equiv 1(mod 2)$.\\
Let ${w_1,w_2,r} \in {\mathbb N}$ with $w_1\equiv 1 (mod 2)$ and $w_2\equiv 1 (mod 2)$ and $h \in {\mathbb Z}$. First, we observe that\\
\begin{equation}\begin{split}
&\frac{1}{[2]_{q^{w_1}}^r} l_{q^{w_1},r}^{(h)}\left(s, w_2 x+\frac{w_2}{w_1}\sum_{l=1}^r j_l|\chi\right)\\
&=\sum_{m_1,\cdot\cdot\cdot,m_r=0}^\infty \frac{(-1)^{m_1+\cdot\cdot\cdot+m_r}q^{w_1\sum_{l=1}^{r}(h-l+1)m_l}(\prod_{l=1}^r\chi(m_l))}
{[m_1+\cdot\cdot\cdot+m_r+{w_2}x+\frac{w_2}{w_1}(j_1+\cdot\cdot\cdot+j_r)]_{q^{w_1}}^s}\\
&=\sum_{m_1,\cdot\cdot\cdot,m_r=0}^\infty \frac{q^{w_1\sum_{l=1}^{r}(h-l+1)m_l}(-1)^{\sum_{l=1}^r m_l}(\prod_{l=1}^r\chi(m_l))[w_1]_q^s}{[w_2(j_1+\cdot\cdot\cdot+j_r)+{w_1}{w_2}x+{w_1}(m_1+\cdot\cdot\cdot+m_r)]_q^s}\\
 &=[w_1]_q^s \sum_{n_1,\cdot\cdot\cdot,n_r=0}^\infty \sum_{i_1,\cdot\cdot\cdot,i_r=0}^{d{w_2}-1}   \frac{(-1)^{\sum_{l=1}^r {(i_l+n_l)}} q^{{w_1}{\sum_{l=1}^r {(h-l+1)(i_l+{n_l}db)}}} (\prod_{l=1}^r \chi(i_l))} {[{w_1}{w_2}(x+d{\sum_{l=1}^r {n_l}})+{w_2}{\sum_{l=1}^r {j_l}}+{w_1} {\sum_{l=1}^r {i_l}}]_q^s}\\
 &=[w_1]_q^s \sum_{n_1,\cdot\cdot\cdot,n_r=0}^\infty \sum_{i_1,\cdot\cdot\cdot,i_r=0}^{{w_2}d-1}   \frac{(-1)^{\sum_{l=1}^r {(i_l+n_l)}} q^{{w_1}{\sum_{l=1}^r {(h-l+1)(i_l+{n_l}{w_2}d)}}} (\prod_{l=1}^r \chi(i_l))} {[{w_1}{w_2}(x+d{\sum_{l=1}^r {n_l}})+{w_2}{\sum_{l=1}^r {j_l}}+{w_1} {\sum_{l=1}^r {i_l}}]_q^s}.\\
\end{split}\end{equation}
Thus, by (2.1), we get
\begin{equation}\begin{split}
&\frac{[w_2]_q^s}{[2]_{q^{w_1}}^r} \sum_{i_1,\cdot\cdot\cdot,i_r=0}^{d{w_1}-1} (-1)^{\sum_{l=1}^r j_l} \left(\prod_{l=1}^r \chi(j_l)\right) q^{b\sum_{l=1}^r (h-l+1){j_l}} l_{q^{w_1},r}^{(h)}\left(s, w_2 x+\frac{w_2}{w_1} \sum_{l=1}^r j_l|\chi\right)\\
&=[w_1]_q^s [w_2]_q^s \sum_{i_1,\cdot\cdot\cdot,i_r=0}^{d{w_2}-1} \sum_{j_1,\cdot\cdot\cdot,j_r=0}^{d{w_1}-1} \sum_{n_1,\cdot\cdot\cdot,n_r=0}^\infty \\
&\left( \frac{(-1)^{\sum_{l=1}^r (i_l+j_l+n_l)} (\prod_{l=1}^r \chi(j_l))(\prod_{l=1}^r \chi(i_l)) q^{w_2\sum_{l=1}^{r}(h-l+1)j_l +w_1\sum_{l=1}^{r}(h-l+1)i_l}}  {[{w_1}{w_2}(x+d{\sum_{l=1}^r {n_l}})+{w_2}{\sum_{l=1}^r {j_l}}+{w_1} {\sum_{l=1}^r {i_l}}]_q^s}  \right)\\
&\times q^{{w_1}{w_2}d{\sum_{l=1}^r (h-l+1){n_l}}}.\\
\end{split}\end{equation}
By the same method as (2.2), we get
\begin{equation}\begin{split}
&\frac{[w_1]_q^s}{[2]_{q^{w_2}}^s} \sum_{j_1,\cdot\cdot\cdot,j_r=0}^{d{w_2}-1} (-1)^{\sum_{l=1}^r j_l} \left(\prod_{l=1}^r \chi(j_l)\right) q^{{w_1}\sum_{l=1}^r (h-l+1){j_l}} l_{q^{w_2},r}^{(h)} \left(s, w_1 x+\frac{w_1}{w_2} \sum_{l=1}^r j_l|\chi\right)\\
&=[w_1]_q^s [w_2]_q^s \sum_{j_1,\cdot\cdot\cdot,j_r=0}^{d{w_2}-1} \sum_{i_1,\cdot\cdot\cdot,i_r=0}^{d{w_1}-1} \sum_{n_1,\cdot\cdot\cdot,n_r=0}^\infty\\
  &\left(\frac{(-1)^{\sum_{l=1}^r (i_l+j_l+n_l)} (\prod_{l=1}^r \chi(j_l))(\prod_{l=1}^r \chi(i_l)) q^{w_1\sum_{l=1}^{r}(h-l+1)j_l +w_2\sum_{l=1}^{r}(h-l+1)i_l}}  {[{w_1}{w_2}(x+d{\sum_{l=1}^r {n_l}})+{w_1}{\sum_{l=1}^r {j_l}}+{w_2} {\sum_{l=1}^r {i_l}}]_q^s}\right)\\
  &\times q^{{w_1}{w_2}d{\sum_{l=1}^r (h-l+1){n_l}}}.\\
\end{split}\end{equation}

Therefore, by(2.2) and (2.3), we obtain the following Theorem.
\bigskip
\begin{thm}\label{Theorem 1.} For ${w_1,w_2} \in{\mathbb N}$ with $w_1\equiv1 (mod2)$ and $w_2\equiv1 (mod2)$, we have\\
\begin{equation*}\begin{split}
&\ [2]_{q^{w_2}}^r [w_2]_q^s \sum_{j_1,\cdot\cdot\cdot,j_r=0}^{{w_1}d-1} (-1)^{\sum_{l=1}^r j_l}\left(\prod_{l=1}^r \chi(j_l)\right) q^{{w_2}\sum_{l=1}^r (h-l+1){j_l}} l_{q^{w_1},r}^{(h)}\left(s, w_2 x+\frac{w_2}{w_1} \sum_{l=1}^r j_l|\chi\right)\\
=&[2]_{q^{w_1}}^r [w_1]_q^s \sum_{j_1,\cdot\cdot\cdot,j_r=0}^{{w_2}d-1} (-1)^{\sum_{l=1}^r j_l}\left(\prod_{l=1}^r \chi(j_l)\right) q^{{w_1}\sum_{l=1}^r (h-l+1){j_l}} l_{q^{w_2},r}^{(h)}\left(s, w_1 x+\frac{w_1}{w_2} \sum_{l=1}^r j_l|\chi\right).\\
\end{split}\end{equation*}
\end{thm}
\bigskip
By(1.7) and Theorem $2.1$, we obtain the following theorem.
\bigskip
\begin{thm}\label{Theorem 2.}  For $n \in{\mathbb Z_{\geq0}}$ and ${w_1,w_2}\in{\mathbb N}$ with $w_1\equiv1 (mod2)$ and $w_2\equiv1 (mod2)$, we have\\
\begin{equation*}\begin{split}
&\ [2]_{q^{w_2}}^r [w_1]_q^n \sum_{j_1,\cdot\cdot\cdot,j_r=0}^{{w_1}d-1} (-1)^{\sum_{l=1}^r j_l}\left(\prod_{l=1}^r \chi(j_l)\right) q^{{w_2}\sum_{l=1}^r (h-l+1){j_l}}  E_{n,\chi,q^{w_1}}^{(h,r)} \left({w_2}x+ \frac{w_2}{w_1} \sum_{l=1}^r j_l\right)\\
=&[2]_{q^{w_1}}^r [w_2]_q^n \sum_{j_1,\cdot\cdot\cdot,j_r=0}^{{w_2}d-1} (-1)^{\sum_{l=1}^r j_l}\left(\prod_{l=1}^r \chi(j_l)\right) q^{{w_1}\sum_{l=1}^r (h-l+1){j_l}}  E_{n,\chi,q^{w_2}}^{(h,r)} \left({w_1}x+ \frac{w_1}{w_2} \sum_{l=1}^r j_l\right).\\
\end{split}\end{equation*}
\end{thm}
\bigskip
From(1.5), we note that\\
\begin{equation}\begin{split}
E_{n,\chi,q}^{(h,r)}(x+y)&=(q^{x+y} E_{\chi,q}^{(h,r)} +[x+y]_q)^n\\
&=(q^{x+y} E_{\chi,q}^{(h,r)}+q^{x}[y]_q+[x]_q)^n\\
&=\sum_{i=0}^{n} {{n}\choose {i}} q^{xi} (q^y E_{\chi,q}^{(h,r)}+[y]_q)^i [x]_q^{n-i}\\
&=\sum_{i=0}^{n} {{n}\choose {i}} q^{xi} E_{i,\chi,q}^{(h,r)}(y) [x]_q^{n-i}.\\
\end{split}\end{equation}
By(2.4), we get\\
\begin{equation}\begin{split}
&\sum_{j_1,\cdot\cdot\cdot,j_r=0}^{d{w_1}-1} (-1)^{\sum_{l=1}^r j_l} q^{{w_2}\sum_{l=1}^r (h-l+1){j_l}} \left(\prod_{l=1}^r \chi(j_l)\right) E_{n,\chi,q^{w_1}}^{(h,r)} \left({w_2}x+\frac{w_2}{w_1} \sum_{l=1}^r j_l\right)\\
=&\sum_{j_1,\cdot\cdot\cdot,j_r=0}^{d{w_1}-1} (-1)^{\sum_{l=1}^r j_l} q^{{w_2}\sum_{l=1}^r (h-l+1){j_l}} \left(\prod_{l=1}^r \chi(j_l)\right) \sum_{i=0}^{n} {{n}\choose {i}} q^{i{w_2}(j_1+\cdot\cdot\cdot+j_r)}\\
&\times  E_{i,\chi,q^{w_1}}^{(h,r)}({w_2}x)  \left[\frac{{w_2}{(j_1+\cdot\cdot\cdot+j_r)}}{{w_1}}\right]_{q^{w_1}}^{n-i}\\
=&\sum_{j_1,\cdot\cdot\cdot,j_r=0}^{d{w_1}-1} (-1)^{\sum_{l=1}^r j_l} q^{{w_2}\sum_{l=1}^r (h-l+1){j_l}} \left(\prod_{l=1}^r \chi(j_l)\right) \sum_{i=0}^{n} {{n}\choose {i}} q^{(n-i){w_2} {\sum_{l=1}^r j_l}} \\
&\times E_{n-i,\chi,q^{w_1}}^{(h,r)}({w_2}x) \left[\frac{w_2}{w_1} {\sum_{l=1}^r{j_l}}\right]_{q^{w_1}}^{i}\\
=&\sum_{i=0}^{n} {{n}\choose {i}} \left(\frac{[w_2]_q}{{[w_1]}_q}\right)^i E_{n-i,\chi,q^{w_1}}^{(h,r)} ({w_2}x) \sum_{j_1,\cdot\cdot\cdot,j_r=0}^{d{w_1}-1} (-1)^{\sum_{l=1}^r j_l} q^{{w_2} \sum_{l=1}^r (h-l+n-i+1){j_l}}\\
&\times \left(\prod_{l=1}^r \chi(j_l)\right) {\left[j_1+\cdot\cdot\cdot+j_r\right]}_{q^{w_2}}^i\\
=&\sum_{i=0}^{n} {{n}\choose {i}} \left(\frac{[w_2]_q}{[w_1]_q}\right)^i E_{n-i,\chi,q^{w_1}}^{(h,r)} ({w_2}x) S_{n,i,q^{w_2}}^{(h,r)} (w_1d|\chi),
\end{split}\end{equation}
where\\
\begin{equation}\begin{split}
S_{n,i,q}^{(h,r)} (w|\chi) =\sum_{j_1,\cdot\cdot\cdot,j_r=0}^{w-1} (-1)^{\sum_{l=1}^r j_l} q^{\sum_{l=1}^r (h-l+n-i+1){j_l}} [j_1+\cdot\cdot\cdot+j_r]_q^i \left(\prod_{l=1}^r \chi(j_l)\right).\\
\end{split}\end{equation}
From (2.5), we have\\
\begin{equation}\begin{split}
&[2]_{q^{w_2}}^r [w_1]_q^n \sum_{j_1,\cdot\cdot\cdot,j_r=0}^{d{w_1}-1} (-1)^{\sum_{l=1}^r j_l} q^{w_2\sum_{l=1}^r (h-l+1){j_l}} \left(\prod_{l=1}^r \chi(j_l)\right)\\
 &\times E_{n,\chi,q^{w_1}}^{(h,r)} \left({w_2}x+\frac{w_2}{w_1}(j_1+\cdot\cdot\cdot+j_r) \right)\\
=&[2]_{q^{w_2}}^r \sum_{i=0}^{n} {{n}\choose {i}} [w_1]_q^{n-i} [w_2]_q^i E_{n-i,\chi,q^{w_1}}^{(h,r)} ({w_2}x) S_{n,i,q^{w_2}}^{(h,r)} ({w_1}d|\chi).\\
\end{split}\end{equation}
By the same method as (2.7), we get\\
\begin{equation}\begin{split}
&[2]_{q^{w_1}}^r [w_2]_q^n \sum_{j_1,\cdot\cdot\cdot,j_r=0}^{d{w_2}-1} (-1)^{\sum_{l=1}^r j_l} q^{w_1\sum_{l=1}^r (h-l+1){j_l}} \left(\prod_{l=1}^r \chi(j_l)\right)\\
 &\times E_{n,\chi,q^{w_2}}^{(h,r)} \left({w_1}x+\frac{w_1}{w_2}\sum_{l=1}^r j_l\right)\\
=&[2]_{q^{w_1}}^r \sum_{i=0}^{n} {{n}\choose {i}} [w_2]_q^{n-i} [w_1]_q^i E_{n-i,\chi,q^{w_2}}^{(h,r)} ({w_1}x) S_{n,i,q^{w_1}}^{(h,r)} ({w_2}d|\chi).\\
\end{split}\end{equation}
Therefore, by (2.7) and (2.8), we obtain the following theorem.\\
\bigskip
\begin{thm}\label{Theorem 2.}  For $n \in{\mathbb Z_{\geq0}}$ and ${w_1,w_2}\in{\mathbb N,}$ with $w_1\equiv1 (mod2)$ and $w_2\equiv1 (mod2)$, we have\\
\begin{equation*}\begin{split}
&[2]_{q^{w_2}}^r \sum_{i=0}^{n} {{n}\choose {i}} [w_1]_q^{n-i} [w_2]_q^i E_{n-i,\chi,q^{w_1}}^{(h,r)} ({w_2}x) S_{n,i,q^{w_2}}^{(h,r)} ({w_1}d|\chi)\\
=&[2]_{q^{w_1}}^r \sum_{i=0}^{n} {{n}\choose {i}} [w_2]_q^{n-i} [w_1]_q^i E_{n-i,\chi,q^{w_2}}^{(h,r)} ({w_1}x) S_{n,i,q^{w_1}}^{(h,r)} ({w_2}d|\chi).\\
\end{split}\end{equation*}
\end{thm}
\bigskip
Now, we observe that\\
\begin{equation}\begin{split}
& e^{{[x]_q}u} \sum_{m_1,\cdot\cdot\cdot,m_r=0}^\infty  q^{\sum_{l=1}^r (h-l+1){m_l}} (-1)^{\sum_{l=1}^r m_l} \left(\prod_{l=1}^r \chi(m_l)\right) e^{[y+{\sum_{l=1}^r m_l}]_q q^x(u+v)}\\
=& e^{{-[x]_q}u} \sum_{m_1,\cdot\cdot\cdot,m_r=0}^\infty  q^{\sum_{l=1}^r(h-l+1){m_l}} (-1)^{\sum_{l=1}^r m_l} \left(\prod_{l=1}^r \chi(m_l)\right) e^{[x+y+{\sum_{l=1}^r m_l}]_q (u+v)}.\\
\end{split}\end{equation}
The left hand side of (2.9) multiplied by $[2]_q^r$ is given by\\
\begin{equation}\begin{split}
&[2]_q^r e^{{[x]_q}u} \sum_{m_1,\cdot\cdot\cdot,m_r=0}^\infty  q^{\sum_{l=1}^r (h-l+1){m_l}} (-1)^{\sum_{l=1}^r m_l} e^{[y+{\sum_{l=1}^r m_l}]_q q^x(u+v)}\left(\prod_{l=1}^r \chi(m_l)\right)\\
=&e^{{[x]_q}u} \sum_{n=0}^\infty q^{nx} E_{n,\chi,q}^{(h,r)}(y) \frac{(u+v)^n}{n!}\\
=& \left(\sum_{l=0}^\infty [x]_q^l \frac{u^l}{l!}\right)\left(\sum_{k=0}^\infty \sum_{n=0}^\infty q^{(k+n)x} E_{k+n,\chi,q}^{(h,r)}(y) \frac{u^k}{k!} \frac{v^n}{n!}\right)\\
=& \sum_{m=0}^\infty \sum_{n=0}^\infty \left(\sum_{k=0}^{m} {{m}\choose{k}} q^{(k+n)x} E_{k+n,\chi,q}^{(h,r)}(y) [x]_q^{m-k}\right) \frac{u^m}{m!} \frac{v^n}{n!}.\\
\end{split}\end{equation}
The right hand side of (2.9) multiplied by $[2]_q^r$ is given by\\
\begin{equation}\begin{split}
&[2]_q^r e^{{-[x]_q}v} \sum_{m_1,\cdot\cdot\cdot, m_r=0}^\infty (-1)^{\sum_{l=1}^r m_l} q^{\sum_{l=1}^r (h-l+1){m_l}} \left(\prod_{l=1}^r \chi(m_l)\right) e^{[x+{\sum_{l=1}^r m_l}]_q (u+v)}\\
=&e^{{-[x]_q}v} \sum_{n=0}^\infty E_{n,\chi,q}^{(h,r)}(x+y) \frac{(u+v)^n}{n!}\\
=&\left(\sum_{l=0}^\infty \frac{({-[x]_q})^l}{l!} v^l\right)\left(\sum_{m=0}^\infty \sum_{k=0}^\infty E_{m+k,\chi,q}^{(h,r)}(x+y) \frac{u^m}{m!} \frac{v^k}{k!}\right)\\
=&\sum_{n=0}^\infty \sum_{m=0}^\infty \left(\sum_{k=0}^{n} {{n}\choose {k}} E_{m+k,\chi,q}^{(h,r)}(x+y) (-[x]_q)^{n-k}\right) \frac{u^m}{m!} \frac{v^n}{n!}\\
=&\sum_{n=0}^\infty \sum_{m=0}^\infty \left(\sum_{k=0}^{n} {{n}\choose {k}} E_{m+k,\chi,q}^{(h,r)}(x+y) q^{(n-k)x} [-x]_q^{n-k}\right) \frac{u^m}{m!} \frac{v^n}{n!}.\\
\end{split}\end{equation}
Therefore, by(2.10) and (2.11), we obtain the following theorem.\\
\bigskip
\begin{thm}\label{Theorem 3.}  For ${m,n}\geq0$, we have
\begin{equation*}\begin{split}
&\sum_{k=0}^{m} {{m}\choose {k}} q^{kx}  E_{n+k,\chi,q}^{(h,r)}(y) {[x]_q}^{m-k}\\
=&\sum_{k=0}^{n} {{n}\choose {k}} q^{-kx} E_{m+k,\chi,q}^{(h,r)}(x+y) [-x]_q^{n-k}.\\
\end{split}\end{equation*}
\end{thm}

\end{document}